\documentclass[11pt,fullpage]{article} 
\usepackage[centertags]{amsmath}
\usepackage{amsfonts}
\usepackage{amssymb}
\usepackage{amsthm}
\usepackage{eepic}
\usepackage{newlfont}
\usepackage{algorithm}
\usepackage{algorithmic}
\usepackage[inactive]{srcltx}  
\theoremstyle{plain}
\newtheorem{thm}{Theorem}[section]
\newcommand{\BTHM}{\begin{thm}} \newcommand{\ETHM}{\end{thm}}
\newtheorem{cor}[thm]{Corollary}
\newcommand{\BCR}{\begin{cor}} \newcommand{\ECR}{\end{cor}}
\newtheorem{lem}[thm]{Lemma}
\newcommand{\BL}{\begin{lem}}   \newcommand{\EL}{\end{lem}}
\newtheorem{clm}[thm]{Claim}
\newcommand{\BCM}{\begin{clm}}   \newcommand{\ECM}{\end{clm}}
\newtheorem{prop}[thm]{Proposition}
\newcommand{\BP}{\begin{prop}}   \newcommand{\EP}{\end{prop}}
\newtheorem{assm}[thm]{Assumption}
\newcommand{\BASM}{\begin{assm}}   \newcommand{\EASM}{\end{assm}}

\theoremstyle{definition}
\newtheorem{defn}{Definition}[section]
\newcommand{\BD}{\begin{defn}}   \newcommand{\ED}{\end{defn}}
\newtheorem{con}[thm]{Conjecture}
\newcommand{\BCONJ}{\begin{con}}   \newcommand{\ECONJ}{\end{con}}

\theoremstyle{definition}
\newtheorem{problem}[thm]{Problem}
\newcommand{\BPR}{\begin{problem}}   \newcommand{\EPR}{\end{problem}}

\newenvironment{rem}{\noindent{\bf Remark:~~}}{}
\newcommand{\BREM}{\begin{rem}} \newcommand{\EREM}{\end{rem}}
\newenvironment{discussion}{\noindent{\bf Discussion:~~\\}}{}
\newcommand{\BDIS}{\begin{discussion}} \newcommand{\EDIS}{\end{discussion}}

\numberwithin{equation}{section}

\def\blackslug
     {\hbox{\hskip 1pt\vrule width 8pt height 8pt depth 1.5pt\hskip 1pt}}
\def\qed{\quad\blackslug\lower 8.5pt\null\par}

\newtheorem{exmp}[thm]{Example}
\newcommand{\BEX}{\begin{exmp}} \newcommand{\EEX}{\end{exmp}}
\newcommand{\BF}{\begin{fact}}   \newcommand{\EF}{\end{fact}}
\newcommand{\Bcr}{\begin{techcorr}}
\newcommand{\Ecr}{\end{techcorr}}
\newcommand{\BDS}{\begin{description}}
\newcommand{\EDS}{\end{description}}
\newcommand{\BE}{\begin{enumerate}}
\newcommand{\EE}{\end{enumerate}}
\newcommand{\BI}{\begin{itemize}}
\newcommand{\EI}{\end{itemize}}
\renewenvironment{proof}{\noindent{\bf Proof:~~}}{\qed}
\newcommand{\BPF}{\begin{proof}}
\newcommand{\EPF}{\end{proof}}

\newcommand{\BB}{\begin{enumerate}}
\newcommand{\EB}{\end{enumerate}}

\title{A note on a Caro-Wei bound for the bipartite independence number in graphs}

\author{Shimon Kogan \\ \\
  Department of Computer Science and Applied Mathematics \\
          Weizmann Institue, Rehovot 76100, Israel 
          \\
          \ shimon.kogan@weizmann.ac.il
          }

\begin{document}
\maketitle
\begin{abstract}
A bi-hole of size $t$ in a bipartite graph $G$ is a copy of $K_{t,t}$ in the bipartite complement of $G$. Given an $n \times n$  bipartite graph $G$, let $\beta(G)$ be the largest $k$ for which 
$G$ has a bi-hole of size $k$. 
We prove that
\[
\beta(G) \geq \left \lfloor \frac{1}{2} \cdot \sum_{v \in V(G)} \frac{1}{d(v)+1} \right \rfloor.
\]
Furthermore, we prove the following generalization of the result above. 
Given an $n \times n$  bipartite graph $G$, let $\beta_d(G)$ be the largest $k$ for which $G$ has
a $k \times k$ induced $d$-degenerate subgraph. We prove that 
\[
\beta_d(G) \geq \left \lfloor \frac{1}{2} \cdot \sum_{v \in V(G)} \min\left(1,\frac{d+1}{d(v)+1}\right) \right \rfloor.
\]
Notice that $\beta_0(G) = \beta(G)$.

\end{abstract}




\section{Introduction}
A bi-hole of size $t$ in a bipartite graph $G$ is a copy of $K_{t,t}$ in the bipartite complement of $G$. Given an $n \times n$  bipartite graph $G$, let $\beta(G)$ be the largest $k$ for which 
$G$ has a bi-hole of size $k$. Denote by $d(v)$ the degree of vertex $v$ in graph $G$.
 The following theorem is proven in \cite{ehard2020biholes}.
\BTHM\label{bihole2}
Given an  $n \times n$  bipartite graph $G$ with average degree $d$ we have
\[
\beta(G) \geq \frac{n}{d+1} - 2.
\]
\ETHM
We generalize and improve the theorem above by proving the following.
\BTHM\label{bipartitecarowei}
Given an  $n \times n$  bipartite graph $G$ we have
\[
\beta(G) \geq \left \lfloor \frac{1}{2} \cdot \sum_{v \in V(G)} \frac{1}{d(v)+1} \right \rfloor.
\]
\ETHM
Notice that Theorem \ref{bihole2} follows from Theorem \ref{bipartitecarowei} by Jensen's inequality. 
The result above is somewhat similar to the well known Caro-Wei bound
\cite{caro2,wei1}, stated below.
\BTHM
Let $\alpha(G)$ denote the maximum number of vertices of
an independent set of graph $G$. Then given a graph $G$, we have the following inequality.
 \[
\alpha(G) \geq  \sum_{v \in V(G)} \frac{1}{d(v)+1}.
\]
\ETHM
We note that better results than the ones in Theorem \ref{bipartitecarowei} exist when the average degree is large. For example, 
in \cite{DBLP:journals/jgt/FeigeK10} the following theorem was proven.
\BTHM\label{logthm}
For any $0 < \epsilon < 1$ there is a constant $d_0$ such that the following holds.
Given an  $n \times n$  bipartite graph $G$ with average degree $d \geq d_0$ and
 $n \geq (1 + \epsilon)d$ we have
\[
\beta(G) \geq \frac{\epsilon}{2} \cdot \frac{n \ln d}{d}.
\] 
\ETHM
Theorems similar to Theorem \ref{logthm} were proven in \cite{axenovich2020bipartite} and \cite{ehard2020biholes}. \\
Now we shall discuss a generalization of Theorem \ref{bipartitecarowei}.
A graph $H$ is $d$-degenerate if
every non-empty subgraph of it contains a vertex of degree at most $d$. 
Thus $0$-degenerate graphs are independent sets and $1$-degenerate graphs are forests.
Let $\alpha_d(G)$ denote the maximum number of vertices of
an induced $d$-degenerate subgraph of $G$.
The following theorem was proven in \cite{DBLP:journals/gc/AlonKS87}.
\BTHM\label{alondegeneratethm1}
 Let $G$ be a graph. Then 
 \[
\alpha_d(G) \geq  \sum_{v \in V(G)} \min\left(1,\frac{d+1}{d(v)+1}\right).
\]
\ETHM
The special case $d=0$ of Theorem \ref{alondegeneratethm1} is the well known Caro-Wei bound \cite{caro2,wei1}. 
We prove an analogous theorem to Theorem \ref{alondegeneratethm1} in bipartite graphs for biholes.
Given an $n \times n$  bipartite graph $G$, let $\beta_d(G)$ be the largest $k$ for which $G$ has
an induced $k \times k$ $d$-degenerate subgraph. We prove the following.
\BTHM\label{ourdegeneratethm1}
Given an  $n \times n$  bipartite graph $G$ and an integer $d\geq 0$ we have
\[
\beta_d(G) \geq \left \lfloor \frac{1}{2} \cdot \sum_{v \in V(G)} \min\left(1,\frac{d+1}{d(v)+1}\right) \right \rfloor.
\]
\ETHM
Notice that Theorem 
\ref{bipartitecarowei} is the special case of $d=0$ in Theorem \ref{ourdegeneratethm1}. 
We prove Theorem \ref{bipartitecarowei} and Theorem \ref{ourdegeneratethm1} in the next two sections.

\section{ Proof of Theorem  \ref{bipartitecarowei} }
We shall prove the following slightly stronger result. 
Given an  $n \times n$  bipartite graph $G=(A,B,E)$ where $\Delta_A$ is the maximum degree in $G$ of the vertices of $A$ and $\Delta_B$ is the maximum degree in $G$ of the vertices of $B$,  
we have
\[
\beta(G) \geq  \frac{1}{2} \left( \frac{1}{\Delta_A+1} + \frac{1}{\Delta_B+1} + \sum_{v \in A \cup B} \frac{1}{d(v)+1} \right) -  1.
\]
Notice that Theorem \ref{bipartitecarowei} will follow from the fact that
\[
\left\lceil \frac{1}{2} \left( \frac{1}{\Delta_A+1} + \frac{1}{\Delta_B+1} + \sum_{v \in A \cup B} \frac{1}{d(v)+1} \right) -  1 \right\rceil \geq 
\left \lfloor \frac{1}{2} \cdot \sum_{v \in A \cup B} \frac{1}{d(v)+1} \right \rfloor
\]
As $\lceil x-\xi\rceil \ge \lfloor x\rfloor$ for any $\xi\in [0,1)$, and in our case $0\leq \xi = 1 - \frac{1}{2} \left( \frac{1}{\Delta_A+1} + \frac{1}{\Delta_B+1} \right)<1$ and $x=\frac{1}{2} \cdot \sum_{v \in A \cup B} \frac{1}{d(v)+1}$.
 \\
Define the potential function $f(d) = \frac{1}{d+1}$ and let the degree sequence of graph $G$ be 
$d_1,d_2,\ldots,d_{2n}$, 
hence we need to prove that $\beta(G) \geq S$ where
\[
S =  \frac{1}{2} \left( f(\Delta_A) + f(\Delta_B) + \sum_{i=1}^{2n}    f(d_i) \right) -  1.
\]
We prove this claim by induction on $n$.
The base case  $n=1$ is trivially correct. Furthermore we may assume that $\Delta_A \geq 1$ and $\Delta_B \geq 1$ for otherwise the graph is an independent set and the claim follows once again trivially. \\
Now we shall consider two cases. \\
\textbf{Case 1:} There is a vertex $a \in A$ such that $d(a) = \Delta_A$ and a vertex 
$b \in B$ such that $d(b) = \Delta_B$ and there is no edge between vertices $a$ and $b$. \\
Let graph $H(A',B',E')$ be the $(n-1) \times (n-1)$ bipartite graph formed from $G$ by removing vertices $a$ and $b$, and let $d'_1,d'_2,\ldots,d'_{2n-2}$ be the degree sequence of graph $H$.
Let
\[
Q =  \frac{1}{2} \left( f(\Delta_{A'}) + f(\Delta_{B'}) + \sum_{i=1}^{2n-2}    f(d'_i) \right) -  1,
\]
where $\Delta_{A'}$ is the maximum degree in $H$ of the vertices of $A'$ and $\Delta_{B'}$ is the maximum degree in $H$ of the vertices of $B'$. 
Now notice that 
\begin{align}
Q &\geq S - \frac{1}{2}\left( f(\Delta_A) + f(\Delta_B) - \Delta_A 
(   f(\Delta_B-1) - f(\Delta_B) ) - \Delta_B
(   f(\Delta_A-1) - f(\Delta_A) ) \right) \notag \\
&= S - \frac{1}{2} \left( \frac{1}{\Delta_A+1} + \frac{1}{\Delta_B+1} 
- \frac{\Delta_A}{\Delta_B(\Delta_B+1)} - \frac{\Delta_B}{\Delta_A(\Delta_A+1)} 
\right) \notag \\ &= S + \frac{1}{2} \frac{(\Delta_A - \Delta_B)^2 (\Delta_A + \Delta_B +1)}{\Delta_A (\Delta_A +1) \Delta_B (\Delta_B+1)} \notag \\ &\geq S \notag
\end{align}
And thus we are done by applying the induction hypothesis to graph $H$. \\
\textbf{Case 2:} Each vertex $a \in A$ such that $d(a) = \Delta_A$ and each vertex 
$b \in B$ such that $d(b) = \Delta_B$ are joined by an edge. \\
Pick an arbitrary vertex  $a \in A$ such that $d(a) = \Delta_A$ and an arbitrary vertex 
$b \in B$ such that $d(b) = \Delta_B$. Note that there is an edge between $a$ and $b$. 
Let graph $H(A',B',E')$ be the $(n-1) \times (n-1)$ bipartite graph formed from $G$ by removing vertices $a$ and $b$, and let $d'_1,d'_2,\ldots,d'_{2n-2}$ be the degree sequence of graph $H$.
Let
\[
Q =  \frac{1}{2} \left( f(\Delta_{A'}) + f(\Delta_{B'}) + \sum_{i=1}^{2n-2}    f(d'_i) \right) -  1,
\]
where $\Delta_{A'}$ is the maximum degree in $H$ of the vertices of $A'$ and $\Delta_{B'}$ is the maximum degree in $H$ of the vertices of $B'$. 
Notice that  by the definition of case $2$ we have $\Delta_{A'} \leq \Delta_A - 1$ and $\Delta_{B'} \leq \Delta_B - 1$. 
Hence we have
\begin{align}
Q &\geq S - \frac{1}{2}\left( f(\Delta_A) + f(\Delta_B) - (\Delta_A-1) 
(   f(\Delta_B-1) - f(\Delta_B) ) - (\Delta_B-1)
(   f(\Delta_A-1) - f(\Delta_A) ) \right) \notag \\
&\qquad + \frac{1}{2}\left( (f(\Delta_{A'}) -  f(\Delta_{A})) + (f(\Delta_{B'}) -  f(\Delta_{B}))  \right) \notag \\
&\geq S - \frac{1}{2} \left( \frac{1}{\Delta_A+1} + \frac{1}{\Delta_B+1} 
- \frac{\Delta_A}{\Delta_B(\Delta_B+1)} - \frac{\Delta_B}{\Delta_A(\Delta_A+1)} 
\right) \notag \\ &= S + \frac{1}{2} \frac{(\Delta_A - \Delta_B)^2 (\Delta_A + \Delta_B +1)}{\Delta_A (\Delta_A +1) \Delta_B (\Delta_B+1)} \notag \\ &\geq S \notag
\end{align}
And thus we are done by applying the induction hypothesis to graph $H$.

\section{ Proof of Theorem  \ref{ourdegeneratethm1} }
We shall prove the following slightly stronger result. Set a fixed integer $d \geq 1$
(the $d=0$ case is Theorem \ref{bipartitecarowei}). 
Define the potential function \[f(x) = \min \left( 1, \frac{d+1}{x+1} \right).\] Let the degree sequence of 
 the $n \times n$ bipartite graph $G=(A,B,E)$ be 
$d_1,d_2,\ldots,d_{2n}$.  Furthermore let $\Delta_A$ be the maximum degree in $G$ of the vertices of $A$ and $\Delta_B$ be the maximum degree in $G$ of the vertices of $B$. We claim
that $\beta_d(G) \geq S$ where
\[
S =  \frac{1}{2} \left( f(\Delta_A) + f(\Delta_B) + \sum_{i=1}^{2n} f(d_i) \right) -  1.
\]
We prove this claim by induction on $n$.
The base case is $n=1$ is trivially correct. 
Now we assume that the claim holds for all bipartite graphs on $(n-1) \times (n-1)$ vertices 
and prove it for bipartite graphs on $n \times n$ vertices. We do this by induction on the number of edges in the $n \times n$ bipartite graph $G$. 
The base case when $G$ has no edges follows trivially from the fact that such a graph is an independent set. Now if $G$ contains a vertex $v$ such that $1 \leq d(v) \leq d$ then we delete the edges incident to $v$ and apply the edge induction hypothesis on the resulting graph $G'$.
Hence $G'$ contains a $k \times k$ $d$-degenerate subgraph $H$ such that $k \geq S$. If $H$ does not contain vertex $v$ then $H$ is also a $k \times k$ $d$-degenerate subgraph of $G$. If $H$ contains vertex $v$ then we add back the edges that are incident to $v$ in $G$ and are inside $H$. The resulting graph $H'$ is a  $d$-degenerate subgraph of $G$ since we took a vertex $v$ of degree $0$ in a $d$-degenerate subgraph and added to it at most $d$ edges. \\
Hence we can assume that the minimum degree of each vertex which is not of degree $0$ in  $G$ is at least $d+1$ and furthermore $\Delta_A \geq d+1$ and 
$\Delta_B \geq d+1$. The rest of the proof is identical to the proof of Theorem \ref{bipartitecarowei} and thus omitted (in particular we do the same case analysis of the two cases in the proof of Theorem \ref{bipartitecarowei} but with potential function $f(x) = \min \left( 1, \frac{d+1}{x+1} \right)$).

\section{Concluding remarks}
It would be interesting to improve Theorem \ref{bipartitecarowei} for bipartite graphs without cycles of length $4$. Furthermore, it would be interesting to generalize the results for $k$-partite graphs where $k\geq 3$.

\section*{Acknowledgements}
The work of the author was supported in part by the Israel Science Foundation (grant No. 1388/16).  The author thanks Uriel Feige and Yair Caro for helpful discussions.

\bibliographystyle{alpha}





\end{document}